\documentclass[12pt]{article} 
\usepackage{setspace}
\usepackage[final]{graphicx}
\usepackage[dvips]{color}
\usepackage{amsmath}
\usepackage{amssymb}
\usepackage{ascmac}
\usepackage{amsthm}
\usepackage{bm}
\usepackage{url}
\usepackage{tabularx}
\usepackage{fancyhdr}
\usepackage[dvipdfmx,svgnames]{xcolor}
\usepackage{tikz}
\usepackage{here}
\usetikzlibrary{calc}
\usepackage{pgfplots}
\usetikzlibrary{positioning}
\newcommand\N{\mathbb{N}}
\newcommand\C{\mathbb{C}}

\newcommand\Z{\mathbb{Z}}

\newcommand\al{{\alpha}}

\newcommand\g{{\gamma}}  
\newcommand\G{{\Gamma}} 

\newmuskip\pFqskip
\mathchardef\pFcomma=\mathcode`, 

\newcommand*\pFq[5]{%
  \begingroup
  \begingroup\lccode`~=`,
    \lowercase{\endgroup\def~}{\pFcomma\mkern\pFqskip}%
  \mathcode`,=\string"8000
  {}_{#1}F_{#2}\biggl(\genfrac..{0pt}{}{#3}{#4};#5\biggr)%
  \endgroup
}
\newtheorem{lmm}{Lemma}
\newtheorem{thm}[lmm]{Theorem}

\newtheorem{rmk}[lmm]{Remark}

\def\comment#1{ }

\makeatletter
    
    \@addtoreset{equation}{section}
  \makeatother
\allowdisplaybreaks
\begin{document}
\title{
On a strange evaluation \\of the hypergeometric series by Gosper. II
}
\author{Akihito Ebisu
}

\maketitle
\begin{abstract}
There are many identities for the hypergeometric series
presented in the article ``Special values of the hypergeometric series'' by Ebisu.
In this note,
we obtain a new hypergeometric identity,
which includes some of these identities as special cases.
We notice that
this identity closely relates to a strange evaluation by Gosper.


\end{abstract}
\section{Introduction and Main Theorem}
We begin with the binomial theorem,
\begin{gather*}
\pFq{1}{0}{a}{-}{x}:=\sum _{n=0} ^{\infty} \frac{(a) _n}{n!}x^n
=
(1-x)^{-a},
\label{binom}
\end{gather*}
where $(a) _n$ is the Pochhammer symbol defined by
\begin{gather*}
(a) _n :=\frac{\G (a + n)}{\G (a)}
=
\begin{cases}
 1  & n = 0, \\
 a (a+1)\cdots (a +n-1) & n \in \N.
\end{cases}
\label{poch}
\end{gather*}
Let us consider a generalization of $_1F_0(a ; - ; x)$ defined as
\begin{gather*}
\pFq{2}{1}{a,b}{c}{x}:=
\sum _{n=0} ^{\infty}\frac{(a)_n (b) _n}{(c) _n n!}x^n,
\label{hgs}
\end{gather*}
where $c \notin \Z _{\leq 0}$.
This series is called the hypergeometric series.
Unfortunately, in the general case, with unrestricted values
of $(a, b, c, x)$, the hypergeometric series $_2F_1(a, b; c; x)$
cannot be expressed in terms of well-known functions,
or to be more exact,
gamma functions together with elementary functions
(see \cite{Ze}).
However, 
$_2F_1(a, b; c; x)$ can be evaluated for parameter values
satisfying certain conditions. 
For instance, 
\begin{gather*}
\pFq{2}{1}{a,b}{c}{1}
=
\frac{\G (c) \G (c-a-b)}{\G (c-a) \G (c-b)}
\end{gather*}
holds for $\Re (c-a-b) > 0$.
There are many other known identities for $_2F_1(a, b; c; x)$. 
Most of these identities have been derived using Gosper's algorithm, the W-Z method,
Zeilberger's algorithm
(see \cite{Ko} and \cite{PWZ}),
and the method of contiguity relations, which was recently introduced in \cite{Eb2}.
In \cite{Eb2}, a number of identities for $_2F_1(a, b; c; x)$ are tabulated
(see also \cite{Ek}, \cite{Ge} and \cite{Gs}). 

First, we expand the definition of $_2F_1(a,b;c;x)$.
Even if the parameter $c$ is a non-positive integer,
we define $_2F_1(a,b;c;x)$ as follows 
if $b$ is a non-positive integer satisfying $c < b$:
\begin{gather*}
\pFq{2}{1}{a,b}{c}{x}
:=
\sum _{n=0} ^{|b|}
\frac{(a) _n (b) _n}{(c) _n n!}x^n
\qquad
(b, c \in \Z _{\leq 0};\ c < b).
\label{ghs2}
\end{gather*}
With this expanded definition, 
for example, the following identities hold for any $m \in \Z _{\ge 0}$:
\begin{flalign}
&\pFq{2}{1}{a,3\, a+1}{3\, a}{\frac{3}{2}}
=\begin{cases}
0
&{\text {if $a=-1-m$}},\\
\vspace{2pt}
{\dfrac { \left( -3 \right) ^{3\,m} \left( 1/3 \right) _m
 \left( 5/3 \right) _{2\,m} }{{2}^{3\,m} \left( 2 \right) _{3\,m} }}
&{\text {if $a=-1/3-m$}},\\
{\dfrac { \left( -3 \right) ^{3\,m} \left( 2/3 \right) _{m} \left( 7/3 \right) _{2\,m} }{{2}^{3\,m+1} \left( 3\right) _{3\,m } }}
&{\text {if $a=-2/3-m$}},
\end{cases}&
\label{case1}
\\
&\pFq{2}{1}{a,4\,a+1}{4\,a}{\frac{4}{3}}
=\begin{cases}
0
& {\text {if $a=-1-m$}},\\
\vspace{2pt}
{\dfrac { \left( -1 \right) ^{m}{2}^{8\,m} 
\left( 1/4\right) _{m}  \left( 7/4 \right) _{3\,m} }{{3}^{4\,m} \left( 2 \right) _{4\,m} }}
& {\text {if $a=-1/4-m$}},\\
\vspace{2pt}
{\dfrac { \left( -1 \right) ^{m}{2}^{8\,m+1} 
\left( 1/2 \right) _{m}  \left( 5/2 \right) _{3\,m} }
{{3}^{4\,m+1} \left( 3 \right) _{4\,m} }}
& {\text {if $a=-1/2-m$}},\\
{\dfrac { 5\left( -1 \right) ^{m}{2}^{8\,m-1} 
\left( 3/4 \right) _{m}  \left( 13/4\right) _{3\,m } }
{{3}^{4\,m+2} \left( 4 \right) _{4\,m} }}
 & {\text {if $a=-3/4-m$}}.
\end{cases}&
\label{case2}
\end{flalign}
The formulae (\ref{case1}) and (\ref{case2}) appear 
as (1,3,3-3)(i) and (1,4,4-1)(i), respectively, in \cite{Eb2}. 
These formulae are treated individually in \cite{Eb2}.
However, looking closely at them,
we realize that they have a similarity:
Their left-hand sides both have the form
\begin{gather}
\pFq{2}{1}{\al, 1-k}{-k}{\frac{k}{\al +k}},
\label{gosper1}
\end{gather}
where $k\in \N$.
Hence, if we are able to evaluate (\ref{gosper1}),
then (\ref{case1}) and (\ref{case2}) 
follow as special cases. 

From Gosper's algorithm, we find that
\begin{gather}
  \frac{(\al)_n(1-k)_{n}}{(-k)_nn!}
  \left(
  \frac{k}{\al+k}
  \right)^n
  =
  f(n+1)-f(n)
  \label{algorithm}
\end{gather}
where
\begin{gather*}
  f(n):=
  \frac{\al+k}{k}
  \cdot
  \frac{(\al+1)_{n-1}}{(n-1)!}
  \left(
  \frac{k}{\al+k}
  \right)^n.
\end{gather*}
Formula (\ref{algorithm}) implies 
\begin{gather}
\pFq{2}{1}{\al, 1-k}{-k}{\frac{k}{\al +k}}
=
\frac{(\al +1) _k}{k!}\left(\frac{k}{\al +k}\right)^k,
\label{gosper2}
\end{gather}
and, from this, we have the following theorem:
\begin{thm}
For any $(\al, k)$ satisfying
\begin{gather*}
(\al, k)\in
\{
(\al ,k) \in \C \times \N \ | \  \al +k \neq 0
 \},
\end{gather*}
we have (\ref{gosper2}).
The identity (\ref{gosper2}) includes 
the formulae (\ref{case1}) and (\ref{case2}),
and also the formulae (1,5,5-1)(i),
(2,5,5-1)(i) and
(1,6,6-1)(i) in \cite{Eb2},
as special cases.
\label{thm}
\end{thm}
The above considerations illustrate that
by tabulating and closely inspecting formulae for mathematical functions,
we can sometimes obtain new formulae.

In this note,
we give another proof of the above theorem;
Formula (\ref{gosper2}) can be easily derived by hand,
and we realize that
(\ref{gosper2}) 
closely relates to a strange evaluation
by Gosper
(see formula (\ref{lmm3})).

\section{A Proof of Theorem 1}
In this section, we prove Theorem 1.

The following identity is easily verified:
\begin{gather}
\pFq{2}{1}{\al, 1+\g}{\g}{x}
=
\frac{(\al x-\g x+\g)(1-x)^{-\al-1}}{\g}.
\label{lmm1}
\end{gather}
Now, we consider the case in which $\g$ approaches $-k$, where $k \in \N$.
Then, the left-hand side of (\ref{lmm1}) becomes
\begin{align}
\begin{split}
&\sum _{n=0} ^{k-1}
\frac{(\al) _n (1-k) _n}{(-k) _n n!}x^n\\
&\quad -\frac{(\al)_{k+1}}{k (k+1)!}x^{k+1}
\left[
1+\frac{(\al+k+1)(2) }{(k+2)\cdot 1!}x
+\frac{(\al+k+1) _2  (2) _2}{(k+2) _2\cdot 2!}x^2
+\cdots
\right]\\
&=
\pFq{2}{1}{\al, 1-k}{-k}{x}
-\frac{(\al)_{k+1}}{k (k+1)!}x^{k+1}
\pFq{2}{1}{\al+k+1, 2}{k+2}{x}.
\end{split}
\label{lmm2}
\end{align}

Next, we evaluate
\begin{gather*}
\pFq{2}{1}{\al+k+1, 2}{k+2}{\frac{k}{\al+k}}.
\end{gather*}
From (40) in Section 2.8 of \cite{Erd},
we have
\begin{align}
\begin{split}
&[c-2b+(b-a)x]\, \pFq{2}{1}{a,b}{c}{x}\\
&+b(1-x)\, \pFq{2}{1}{a,b+1}{c}{x}
-(c-b)\, \pFq{2}{1}{a,b-1}{c}{x}=0.
\end{split}
\label{3tr1}
\end{align}
Substituting $(a,b,c)=(\al + k +1, 1, k+2)$ into (\ref{3tr1}),
this becomes
\begin{align}
\begin{split}
[k-(\al +k)x]\, \pFq{2}{1}{\al+k+1 ,1}{k + 2}{x}
+(1-x)\, \pFq{2}{1}{\al + k +1,2}{k + 2}{x}
=(k + 1).
\end{split}
\label{3tr2}
\end{align}
This relation holds near $x=0$.
Now, we carry out an analytic continuation
of each side of (\ref{3tr2}) along a curve
starting at $x=0$ and ending at $x=k/(\al +k)$.
In this way, we obtain
\begin{gather}
\pFq{2}{1}{\al+k+1, 2}{k+2}{\frac{k}{\al+k}}
=
\frac{(\al + k) (k + 1) }{\al}.
\label{lmm3}
\end{gather}
This formula was first derived by Gosper in \cite{Gos}
(see also \cite{Eb2}).
\begin{rmk}
As we now show,
(\ref{lmm3}) holds for
any $(\al, k)$ in the set
\begin{gather} 
\left\{(\al, k)\in \C^2\ ; \
\al \neq 0,\
{\al +k} \neq 0,\ 
k\neq -2, -3, -4, \ldots
\right\}.
\label{set}
\end{gather}
First, we define
\begin{gather*}
F(x):=\pFq{2}{1}{\al+k+1, 2}{k+2}{x}.
\end{gather*}
Then, because the radius of convergence of a non-terminating $F(x)$ is 1,
and because $F(x)$ is a multivalued function,
we cannot uniquely determine its value for $x \in  \C$ satisfying $|x| \geq 1$,
in general.
In other words, the value of $F(x)$ at such values of $x$ is ill-defined.
However, as seen in the above, the value of $F(k/(\al+k))$ is uniquely determined
and, hence, well-defined even if $|k/(\al+k)| \geq 1$.
\end{rmk}
From (\ref{lmm1}), (\ref{lmm2}) and (\ref{lmm3}), we find that 
\begin{align*}
&\lim _{x\rightarrow k/( \al + k)}\lim _{\g \rightarrow -k}
{\text{(lhs of (\ref{lmm1}))}}
=
\pFq{2}{1}{\al, 1-k}{-k}{\frac{k}{\al +k}}
-
\frac{(\al +1) _k}{k!}\left(\frac{k}{\al +k}\right)^k,\\
&\lim _{x\rightarrow k/( \al + k)}\lim _{\g \rightarrow -k}
{\text{(rhs of (\ref{lmm1}))}}=0.
\end{align*}
Thus, 
we have verified Theorem {\ref{thm}} except in the case $\al = 0$.
However, it is obvious that  (\ref{gosper2}) holds when $\al = 0$.
This completes the proof of Theorem {\ref{thm}}.  

\section*{Acknowledgement}
This work is supported by a Grant-in-Aid for JSPS Fellows, JSPS No. 15J00201.

\medskip
\begin{flushleft}
Akihito Ebisu\\
Department of Mathematics\\
Hokkaido University\\
Kita 10, Nishi 8, Kita-ku, Sapporo, 060-0810\\
Japan\\
a-ebisu@math.sci.hokudai.ac.jp
\end{flushleft}

\end{document}